\input amstex

\def\BC{\Bbb C }
\def\BG{\Bbb G }

\def\hd{, \hdots ,}
\def\inv{{}^{-1}}

\def\na{n+a}

\def\ot{\!\otimes\!}

\def\pp#1{\Bbb P^{#1}}
\def\ppp{\Bbb P}

\def\tdim{\text{dim}\,}

\documentstyle{amsppt}
\magnification = 1100
\hsize =15truecm
\hcorrection{.5truein}
\baselineskip =18truept
 \vsize =22truecm
\NoBlackBoxes
\topmatter
\title On the structure of varieties
with degenerate Gauss mappings\endtitle
 \rightheadtext{Varieties
with degenerate Gauss mappings}

 \author
  M.A. Akivis, V.V. Goldberg, and J.M. Landsberg
\endauthor

 \address  
Department of Mathematics,
Jerusalem College of Technology
21 Havaad Haleumi St., P.O.B. 16031,
Jerusalem 91160, Israel\endaddress

\email {akivis\@avoda.jct.ac.il}\endemail

\address Department of Mathematical Sciences,
New Jersey Institute of Technology,
University Heights,  Newark, NJ 07102, U.S.A.\endaddress

\email {vlgold\@m.njit.edu}\endemail
 
\address Laboratoire de Math\'ematiques,
  Universit\'e Paul Sabatier, UFR-MIG
  31062 Toulouse Cedex,
  FRANCE
\endaddress
\email {jml\@picard.ups-tlse.fr }\endemail

\endkeywords
 
\endtopmatter

\document

Let $V=\BC^{N+1}$ and let $X^n\subset\ppp V$ be a variety. Let $x\in X$ be
a smooth point, and let $\tilde T_xX\subset\ppp V$ denote the embedded
tangent projective space to $X$ at $x$. Let 
$$
\align
\gamma : X& \dasharrow\BG (n, \ppp V)\\
x&\mapsto \tilde T_xX\endalign
$$
denote the {\it  Gauss map} of $X$, where $\BG (n, \ppp V)$ denotes the
Grasmannian of $\pp n$'s in $\ppp V$.

In [GH], Griffiths and Harris present a structure theorem for varieties
with {\it degenerate} Gauss mappings,
that is $X$ such that $\tdim \gamma (X)<\tdim X$. Namely, 
such varieties are
\lq\lq built up from cones and developable varieties\rq\rq\  [GH, p. 392].
By \lq\lq built up from\rq\rq\ they appear to mean \lq\lq foliated
by\rq\rq\  and by \lq\lq developable varieties\rq\rq\ they appear
to mean the osculating varieties to a curve.   With these interpretations,
their result appears to be complete for varieties whose Gauss maps
have one-dimensional fibers. For varieties with higher dimensional
fibers, one could generalize \lq\lq built up from\rq\rq\ to
mean either foliated by, or iteratively constructed from, or
some combination of these two and generalize the osculating
varieties of curves to osculating varieties of arbitrary varieties.
Even with this interpretation however,  their result is still incomplete as
general hypersurfaces with degenerate Gauss maps having fibers of dimension
greater than one cannot be built out of cones and osculating varieties. 
In this note we present examples of varieties with degenerate Gauss mappings.
Some of these examples illustrate Girffiths-Harris' structure theorem, and
some 
(see, for example, IIB.) show its incompleteness.

\medpagebreak

  Fixing
$X^n\subset\ppp V$, let $r$ denote the rank of $\gamma$
and set $f=n-r$, the dimension of a general fiber. If
$x\in X$ is a smooth point, we let $F=\gamma\inv \gamma (x)$
denote the fiber of $\gamma$ (which is a $\pp{f}$).
Let $Z_F\subseteq  F\cap X_{sing}$
denote the {\it focus} of $F$, the points where the
image of a desingularization of $X$
is not immersive. $Z_F$ is a codimension
one subset of
$F$ of degree $n-f$. The number of
components of $Z_F$ and the dimension of the varieties
each of these components sweeps out as one varies $F$ furnish
invariants of $X$. 

\bigpagebreak

Here are some examples
of varieties with degenerate Gauss mappings (which are not mutually
exclusive):

\subheading{I. Joins}

Form the join of $k$ varieties
$Y_1\hd Y_k\subset\ppp V$, 
$$
X= S(Y_1\hd Y_k)=\overline{\cup_{y_j\in Y_j }\ppp_{y_1\hd y_k}}  
$$
where $\ppp_{y_1\hd y_k}$ denotes the projective
space spanned by $y_1\hd y_k$ (generically a $\pp{k-1}$).
Note that $\tdim X\leq \Sigma_j \tdim Y_j + (p-1)$
with equality expected.

 Joins have degenerate
Gauss maps with at least $(k-1)$-dimensional fibers
because
  Terracini's lemma (see [Z, II.1.10]) implies that 
the tangent space to $S(Y_1\hd Y_k)$ is constant
along each $ \ppp^{k-1}_{y_1\hd y_k}$.

Two special cases of this construction:  1.  Let  $ L$ be a linear
space. Then $S(Y,L)$ is   a cone over $Y$ with vertex $L$. 
 2.   $Y_j=Y$
for all $j$. Then $X$ is the union of   the secant $\pp{k-1}$'s to $Y$. 

\smallpagebreak

Joins are built out of cones in the sense that one can use e.g.,   the family
of cones over $Y_2$ with vertices   the points of $Y_1$ to sweep out $X$.

\bigpagebreak

\subheading{II. Varieties built from tangent lines}
 
\subheading{IIA. Tangential varieties} Let $Y\subset\ppp V$
be a   variety and let $\tau (Y)\subset\ppp V$
denote the   union of tangent stars
to $Y$. (If $Y$ is smooth, $\tau (Y)$ is the
union of embedded tangent lines to $Y$.)
$\tau (Y)$ has a degenerate Gauss
map with at least one dimensional fibers.    (see [L] for definitions).
One can also take higher osculating varieties of $Y$ which will also have
degenerate Gauss mappings.
Examples IIB and IIC below generalize $\tau (Y)$.

\subheading{IIB. Hyperbands} 
Let $Y\subset\ppp V$ be a smooth variety 
of dimension $m$ and fix $ k= N-m-1$.
For each $y\in Y$, let $L_y\subset \BG(m+k,\ppp V)$
be such that $\tilde T_yY\subset L_y$ and let
$X=\cup_{y\in Y}L_y$. Then $\tdim X\leq N-1$ (with equality
occurring generically) and $X$ will  have
degenerate Gauss map with at least one-dimensional fibers.
Such a variety $X$ is called a hyperband (see [AG, p. 255]).

The hyperbands with fibers of dimension greater than one 
are not built by families of cones 
and developable varieties. So, they are not covered 
by the Griffiths-Harris structure theorem.

 \medpagebreak

One could seek to generalize   tangential varieties in
a different way, namely by taking a subspace of the tangent
lines through each point of $Y$. 
If $x\in \ppp V$ and $v\in T_x\ppp V$,
we let $\pp 1_{x,v}$ denote the line passing through
$x$ with tangent space spanned by $v$.
Let  $\Delta\subset TY$
be a distribution. One could consider the variety
 $X=\cup_{y\in Y, v\in \Delta_y}\pp 1_{y,v}$ consisting of
the union of tangent lines tangent to $\Delta$. In general
$X$ will {\it not} have a degenerate Gauss map, but it will in
some special cases. The case where $Y\subset Z$ and
one takes $X= \cup_{y\in Y}\tilde T_yZ$ is one special
case.
 Here is another construction:

\subheading{IIC. Unions of conjugate spaces}

Let $II=II_{Y,y}\in S^2T^*_yY\ot N_yY$ denote the
projective second fundamental form of $Y$ at $y$ (see [AG], [GH] or
[L] for a definition).

 Let $Y^{n-1}\subset\pp{n+1}$ be
a variety such that at general points there exist
$n-1$ simultaneous eigen-directions for 
the quadrics in its second fundamental form.
This condition holds for
  generic varieties of codimension two.  
(To make the notion of eigen-direction
  precise, choose a nondegenerate quadric in $II$ to identify
$T$ with $T^*$ and consider the quadrics as endomorphisms of $T$.
The result is independent of the choices.)
Let $X^{n}\subset\pp{n+1}$ be the union of one of these families
of embedded tangent lines.  

The directions
     indicated above are called {\it conjugate directions} on
   $Y^{n-1}$.
 
In higher dimensions it is still possible to have a conjugate
direction or conjugate space, but in this case $Y$ must satisfy
a certain exterior differential system. As is shown in [AG, p. 85]  
  local solutions to this system exist and depend on $n(n-1)$
arbitrary functions
of two variables.  

\smallpagebreak

In this case $X$ is the union of the tangential varieties
of the integral curves for the distribution defined by
the conjugate directions to $Y$.

\bigpagebreak

 \subheading{III. Varieties with $f=1$}

 \subheading{IIIA.  Generic varieties with $f=1$}
We say a variety
$X\subset\ppp V$ with $f=1$ is {\it generic among varieties with $f=1$} if
$Z_{F}$ consists of
$n-1$ distinct points and the variety each
point sweeps out is $(n-1)$-dimensional.
The following theorem follows from results in [AG]:

\proclaim{Theorem}  The varieties $X^n\subset\pp\na$
generic among varieties
with $f=1$ are the union of conjugate lines to some variety
$Y^{n-1}$, with a finite number of lines tangent to a general point of $Y$.
\endproclaim
 
\subheading{IIIB. Classification of  $X^3\subset\pp 4$ with $f=1$}
Here $F$ is a $\pp 1$ and
the focus $Z_F$ is of degree two. There are two classes:

\smallpagebreak

\noindent Class 1: $Z_F$ consists of two distinct points, $z_1,z_2$.

1a. (Generic case) Each $z_j$ traces out a surface $S_j$. 
Here $X$ is the dual variety of a $II$-generic surface in 
$\pp{4*}$ (its Gauss image). Locally $X$ may be described as the
union of a family of lines tangent to  conjugate directions on either
surface (one must take the family that corresponds to conjugate directions
on the other surface).  It may be the case that globally
$S_1=S_2$ and there is a unique construction. 

1b. $z_1$ traces out a surface $S$  and $z_2$ traces out a curve
$C$. Here $X$ may be described as the union of a family of
conjugate lines to $S$, where the conjugate lines intersect along
a curve.

1c. Both $z_j$'s trace out curves, $C_j$. In this case
$X=S(C_1,C_2)$.

\smallpagebreak

\noindent Class 2: $Z_F$ is a single point $z$ of multiplicity
two.

2a. $z$ traces out a surface $S$. In this case $S$ will
have a family of asymptotic lines   and
$X$ is the union of the asymptotic lines to $S$.

2b. $z$ traces out a curve $C$.  An example of $X$
in this case is the union of
a family of planes that are tangent to $C$. We
conjecture that this is the only example.

2c. $z$ is fixed, then $X$ is a cone over $z$.

\smallpagebreak

\noindent{\it Acknowledgement}.  We are grateful to J. Piontkowski for
pointing out that our  interpretation of the Griffiths-Harris announcement
in an earlier version of this paper was   too
narrow, and showing us how some of the examples
above fit into the Griffiths-Harris perspective.

\enddocument